\documentclass[12pt]{article}
\usepackage{amssymb,amsmath}
\newtheorem{Theorem}{Theorem}

\newtheorem{Definition}{Definition}

\newtheorem{Lemma}{Lemma}
\newenvironment{proof}{
  \noindent\textbf{Proof}\ }{\hspace*{\fill}
  \begin{math}\Box\end{math}\medskip}

\begin{document}

\title{Some colouring problems for unit-quadrance graphs} \author{Le Anh
Vinh\\
Department of Mathematics\\
Harvard University\\Cambridge, MA 02138
\\
avle@fas.harvard.edu}
\date{\empty}
\maketitle

\begin{abstract}
  The quadrance between two points $A_1 = ( x_1, y_1 )$ and $A_2 = ( x_2, y_2
  )$ is the number $Q ( A_1, A_2 ) = ( x_1 - x_2 )^2 + ( y_1 - y_2 )^2$. Let
  $q$ be an odd prime power and $F_q$ be the finite field with $q$ elements.
  The unit-quadrance graph $D_q$ has the vertex set $F_q^2$, and $X, Y \in
  F_q^2$ are adjacent if and only if $Q ( A_1, A_2 ) = 1$. In this paper, we study some colouring problems for the unit-quadrance graph $D_q$ and discuss some open problems.
\end{abstract}

\section{Introduction}
In \cite{quadrance}, Wildberger introduces a remarkable new approach to trigonometry and Euclidean geometry by replace distance by quadrance and angle by spread, thus allowing the development of Euclidean geometry over any field. The following definition follows from \cite{quadrance}.

\begin{Definition}
  The \textbf{quadrance} $Q ( A_1, A_2 )$ between the points $A_1 = ( x_1,
  y_1)$, and $A_2 = ( x_2, y_2 )$ in $F_q^2$ is the number
  \[ Q ( A_1, A_2 ) := ( x_2 - x_1 )^2 + ( y_2 - y_1 )^2 . \]
\end{Definition}

This approach motives some colouring problems for unit-quadrance graph over any field. Note that in usual $2$-dimensional Euclidean space $R^2$ then the quadrance between $A_1, A_2$ is unit if and only if the distance between $A_1, A_2$ is unit. 

Let $q$ be an odd prime power and $F_q$ be the finite field with $q$ elements. The unit-quadrance graph $D_q$ has the vertex set $F_q^2$, and $X, Y \in F_q^2$ are adjacent if and only if $Q ( A_1, A_2 ) = 1$. This graph (or so-called finite Euclidean graph) was also studied by Medrano et al in \cite{before}. 

This paper aims to outline some results about colouring problems for unit-quadrance graphs and also discuss several interesting open questions.  

\section{The chromatic number}
Let $\chi ( F_q^2 )$ be the chromatic number of graph $D_q$. In  \cite{vinh}, the author showed that

\begin{Theorem}\label{main}
  Suppose that $q = p^n > 3$ where $p$ is an odd prime number then
  \[ 1 + \frac{q - (-1)^{(q-1)/2}}{q^{1/2}} \leqslant \chi ( F_q^2 ) \leqslant
     \frac{p^n + p^{n - 1}}{2} = q ( 1 / 2 + o ( 1 ) ) . \]
\end{Theorem}

For some primes $q$, the graph $D_q$ contains no triangle. We have the following lemma.

\begin{Lemma}\label{no triangle}(\cite{vinh})
Let $q$ be any prime of the form $q = 12k \pm 7$. Then $D_q$ contains no triangle.
\end{Lemma}

\begin{proof}
Suppose that $D_q$ contains a triangle $XYZ$ with $X = ( m, n )$, $Y = ( m + x, n + y )$ and $Z = ( m + x + u, n + y + v )$  for some $m, n, x, y, u, v \in F_q$ then $x^2 + y^2 = u^2 + v^2 = (x+u)^2 + (y+v)^2 = 1$. This implies that $x u + y v = -1/2$.
We have $( x u + y v )^2 + ( x v - y u )^2 = ( x^2 + y^2 ) ( u^2 + v^2)$ so
  \[ ( x v - y u )^2 = 1 - 1/4 = 3/4. \]
For $q = 12 \pm 7$ then $3$ is not square in $F_q$. Hence we have a contradiction. This concludes the proof of the lemma.
\end{proof}

From Theorem \ref{main} and Lemma \ref{no triangle}, if $q \equiv \pm 7$ (mod $12$) is a prime number then $D_q$ is a triangle-free graph with chromatic number $\chi(D_q) \geq q^{1/2}(1+o(1))$. This give us a family of triangle-free graphs with arbitrary high chromatic numbers.

\begin{Theorem} For $q \equiv \pm 7$ (mod $12$) is a prime number then the unit-quadrance graph $D_q$ of order $q^2$ is a triangle-free graph with chromatic number number
\[\chi(D_q) \geq q^{1/2}(1+o(1)).\]
\end{Theorem}

The bound on the independence numbers follows from an estimate of the Lov\'asz $\theta$-function. This function, introduced by Lov\'asz in \cite{lovasz}, can be defined as follows. If $G=(V,E)$ is a graph, an orthonormal labeling of $G$ is a family $(b_v)_{v \in V}$ of unit vectors in an Euclidean space so that if $u$ and $v$ are distinct non-adjacent vertices, then $b_u^tb_v=0$, that is, $b_u$ and $b_v$ are orthogonal. The $\theta$-number $\theta(G)$ is the minimum, over all orthonormal labelings $b_v$ of $G$ and over all unit vectors $c$, of
\[ \max_{v\in V} \frac{1}{(c^tb_v)^2}.\]

It is known that (see \cite{lovasz}) the independence number of $G$ does not exceed $\theta(G)$. Theorem $9$ of \cite{lovasz} asserts that for $d$-regular graphs $G$ with eigenvalues $d= \lambda_1 \geq \lambda_2 \geq \ldots \lambda_n$,
\begin{equation}\label{2} \theta(G) \leq \frac{-n \lambda_n}{\lambda_1 - \lambda_n}. \end{equation}

In \cite{before}, Medrano et al. give a general bound for eigvenvalue of $D_q$.

\begin{Lemma}\label{eigenvalue} (\cite{before}) Let $\lambda \ne \Delta(D_q)$ be any eigenvalue of graph $D_q$ then $| \lambda | \leq q^{1/2}$.
\end{Lemma}

Recall that the graph $D_q$ is a regular graph with degree $\Delta(D_q) = q - (-1)^{(q-1)/2}$ (see \cite{before}). 
From Lemma \ref{eigenvalue} and equation (\ref{2}), we have
\begin{equation} \label{theta}
  \theta ( D_q ) \leqslant \frac{- q^2 ( - q^{1 / 2} )}{q - ( - 1 )^{( q - 1 )
  / 2} - ( - q^{1 / 2} )} \leqslant q^{3 / 2}.
\end{equation}

Thus, the independence number of $D_q$ does not exceed $q^{3/2}$. In other words, $D_q$ contains no $K_3$ and $\bar{D_q}$ contains no $K_{q^{3/2}}$. Therefore, we have a constructive bound for off-diagonal Ramsey number $R(m,3)$. 

\begin{Theorem}
For $q \equiv \pm 7$ (mod $12$) is a prime number then the unit-quadrance graph $D_q$ of order $q^2$ gives a contructive bound for off-diagonal Ramsey number
\[R(m,3) \geq \Omega(m^{2/3}) .\]
\end{Theorem}

This bound matches with the second best bound of Noga Alon in \cite{alon}. The problem of finding better bounds for the chromatic number of graph $D_q$ touches on an important
question in graph theory: what is the greatest possible chromatic number for a
triangle-free regular graph of order n.

\section{Perfect unit-quadrance graph}

A graph $G$ is perfect if and only if for every induced subgraph $H$ of $G$,
$\chi ( H ) = \omega ( H )$ where $\omega ( H )$ as usual denotes the clique
number. We will determine for which $q$ then the unit-quadrance graph $D_q$ is
perfect. We will need some properties of universal geometry.

\begin{Definition}
  A circle in a finite field $F_q$ with center $A_0  = [ x, y ]$ and quadrance $K \in F_q$ is set of all points $X = ( u, v )$ in $F_q \times F_q$ such that
  \[ Q(A_0, X) = K. \]
\end{Definition}
Let $C_i(X)$ denote the circle centered at $X \in F_q^2$ with quadrance $i$. In \cite{vinh3}, the author proved the following result.

\begin{Lemma} \label{khac khong}(\cite{vinh3})
For any $i, j \ne 0$ in $F_q$. Let $X, Y$ be two distinct points in $F_q^2$ such that $k = Q(X,Y)\ne 0 $. Then the number of intersections of two circles $C_i(X)$, $C_j(Y)$ only depends on $i, j$ and $k$. Precisely,  define \[f ( i, j, k ) := i j - \frac{ (i - j - k )^2 }{ 4}.\]
The number of intersection points is $p_{i j}^k$, where
  \begin{equation}\label{51}
  p_{ij}^k =
   \begin{cases}
    0  & \text{if}\; \; f(i,j,k)\; \text{is non-square},\\
    1  & \text{if}\; \; f(i,j,k) = 0,\\
    2 &\text{if}\; \; f(i,j,k)\; \text{is square}.
   \end{cases}
  \end{equation}

\end{Lemma}

\begin{proof}
Suppose that $X = [ m, n ]$ and $Y = [ m + x, n + y ]$ for some $m, n, x, y \in F_q$ then $x^2 + y^2 = k$. Suppose $C_i(X)$ intersects $C_j(Y)$ at $Z$ where $Z = [ m + x + u, n + y + v ]$ for some $u, v \in F_p$. Then we have $u^2 + v^2 = j$ and $( x + u )^2 + ( y + v^2 ) = i$. It implies that
  \[ x u + y v = \frac{i - j - k}{2} . \]
  But we have $( x u + y v )^2 + ( x v - y u )^2 = ( x^2 + y^2 ) ( u^2 + v^2
  )$ so
  \[ ( x v - y u )^2 = k j - \frac{( i - j - k )^2}{4} = i j - \frac{( k - i - j )^2}{4} = f ( i, j, k ) . \]
  If $f ( i, j, k )$ is non-square number in $F_q$ then it is clear that there does not exist
  such $x, y, u, v$, or $p_{i j}^k = 0$. Otherwise, let $\alpha = ( i - j
  - k ) / 2$ and $f ( i, j, k ) = \beta^2$ for $0 \leqslant \beta \leqslant (
  p + 1 ) / 2$ then
  \begin{eqnarray*}
    x v - y u = \pm \beta,\\
    x u + y v = \alpha.
  \end{eqnarray*}
  Solving for $( u, v )$ with respect to $( x, y )$ we have
  \begin{eqnarray*}
    u = ( \alpha x \mp \beta y ) / j,\\
    v = ( \pm \beta y + \alpha x ) / j. 
  \end{eqnarray*}
  If $\beta = 0$ then we have only one $( u, v )$ for each $( x, y )$, but if
  $\beta \neq 0$ then we have two pairs $( u, v )$. It implies (\ref{51}) and concludes the proof of the lemma.
\end{proof}

We are now ready to prove the main result of this section.

\begin{Theorem}
  $D_q$ is not perfect except for $D_3$.
\end{Theorem}

\begin{proof}
  $D_3$ is perfect because $D_3$ is the line graph of $K_{3, 3}$ and line
  graphs of bipartite graphs are perfect.
  
  From Lemma \ref{no triangle}, $D_5$ and $D_7$ have no triangle. So $\omega ( D_5 ) =
  \omega ( D_7 ) = 2$. Since $D_q$ contains an odd cycle of length $q$ for all
  $q$ odd, $\chi ( D_q ) \geqslant 3$ for all $q \geqslant 5$. Thus, we have
  \[ \omega ( D_5 ) \neq \chi ( D_5 ), \omega ( D_7 ) \neq \chi ( D_7 ) . \]
  This implies that $D_5, D_7$ are not perfect.
    
  For $q = 9$, then from Theorem \ref{main}, we have $\chi ( D_9 ) \geqslant 4$. From
  Lemma \ref{khac khong}, $\omega ( D_9 ) \leqslant 3$ (and it is indeed that $\omega (
  D_9 ) = 3$). Thus $\omega ( D_9 ) < \chi ( D_9 )$ or $D_9$ is not perfect.
  
  For $q \geqslant 11$ then from Theorem \ref{main}, we have
  \[ \chi ( D_q ) \geqslant 1 + \frac{q - 1}{\sqrt{q}} > 4, \]
  so $\chi ( D_q ) \geqslant 5$. From Lemma \ref{khac khong}, $\omega ( D_q ) \leqslant
  4$ for all $q$. Thus, $\omega ( D_q ) < \chi ( D_q )$ for $q \geqslant 11$.
  This implies that $D_q$ is not perfect for $q \geqslant 11$.
  This concludes the proof.
\end{proof}

\section{Some open problems}
\subsection{The choice number}

A $k$-list-assignment $L$ to the vertices of a graph $G$ is the assignment of a list, $L(v)$, of at least $k$ colours to every vertex $v$ of $G$. The graph $G$ is $k$-chosable if for every $k$-list-assignment, we can choose a colour for each vertex from its list such that no two adjacent vertices have the same colours. Then the choice number ch$(G)$ of $G$ is the smallest number $k$ such that $G$ is $k$-choosable. For any graph $G$ we have $\left\lceil n/\alpha(G) \right\rceil \leq \chi(G) \leq ch(G)$. We have the following theorem.

\begin{Theorem} Given $\gamma > 0$, there exists a $q_0(\gamma)$ such that for every $q \geq q_0(\gamma)$ we have
\[ \chi(D_q) \leq ch(D_q) \leq \frac{(1+\gamma)q^2}{\log_2q}.  \]
\end{Theorem}

The proof of this theorem is obmited as it is similar to the proof of Theorem 3 in \cite{3}. This theorem gives us a better upper bound for $\chi(D_q)$. The question of finding a better bound or a constructive colouring which matches with this bound is still open.

\subsection{The achromatic number}

The achromatic number $\psi(G)$ of a graph $G$ is the greatest number of colours in a vertex colouring of $G$ such that no two adjacent vertices have the same colour and for any pair of colours, there is at least one edge of $G$ whose endpoints are coloured with this pair of colours. Suppose that we partition $V(D_q)$ into $m$ classes then at least one class has at most $q^2/m$ vertices. Each vertice in this class has degree $\Delta(D_q)$ and there is at least one edge of $G$ whose endpoints are coloured with this pair of colours so
\[ \frac{q^2}{m} \Delta(D_q) \geq m-1.\]
This implies that $\psi(D_q) \leq q^{3/2}$. For the lower bound, it is easy to see that $\psi(D_q) \geq q+1$ (since no two adjacent vertices have the same colour). These bounds however are still far from the truth. 

\subsection{Edge and total colouring}

An edge colouring of a graph $G$ is an assignment of colours to its edges so that no two incident edges have the same colour. The edge-choromatic number $\chi'(G)$ of a graph $G$ is the least number $k$ of colours for which $G$ has an edge colouring with exactly $k$ colours. It is well known that if $G$ is a regular graph of common degree $\Delta$ with odd order then $\chi'(G) = \Delta + 1$. Thus, $\chi'(D_q)=q - (-1)^{(q-1)/2}$. The edge choice number $ch'(G)$ of a graph $G$ is the smallest integer $k$ such that whenever every edge of $G$ is given a list of at least $k$ colours, there exists an edge colouring of $G$ in which every edge receives a colour from its own list and no two incident edges have the same colour. The List Edge Colouring Conjecture states that for every graph $G$, $ch'(G) = \chi'(G)$. The confirmation or rejection of this conjecture for unit-quadrance graph is also an interesting problem.


\begin{thebibliography}{9}


\bibitem{alon}
Alon, N., Explicit Ramsey graphs and orthonormal labelings, \textit{The Electronic Journal of Combinatorics} \textbf{1} (1994), R12, 8pp.

\bibitem{alon1}
Alon, N., Tough Ramsey graphs without short cycles, \textit{Journal of Algebraic Combinatorics} \textbf{4} (1995), 189-195.

\bibitem{3}
Alon, N., Krivelevich, M., Sudakow, B., List colouring of random and pseudo-random graphs, \textit{Combinatorica} \textbf{19} (4) (1999) 453-472.

\bibitem
{hof} Hoffman, A. J., On eigenvalues and colorings of graphs, in \textit{Graph Theoryand Its Applications} (B. Harris, ed.), Acad. Press, 1970, 79-91.

\bibitem
{lovasz}
Lov\'asz, L., On the Shannon capacity of a graph, \textit{IEEE Transactions on Information Theory
IT-25}, (1979), 1-7.

\bibitem
{before} Medrano, A., Myers, P., Stark, H. M., Terras, A., Finite analogues of Euclidean space,
\textit{Journal of Computational and Applied Mathematics}, \textbf{68} (1996), 221-238.


\bibitem
{vinh} Vinh, L. A., On chromatic number of unit-quadrance graphs (or finite Euclidean graph),
submitted to \textit{The Electronic Journal of Combinatorics}.

\bibitem
{vinh1} Vinh, L. A., Quadrance graphs, to appear on
\textit{Australian Mathematical Society Gazette}.

\bibitem
{vinh2} Vinh, L. A., Random walk on hypergroup of circles in a finite field, \textit{The proceedings of Australasian Workshop on Combinatorial Algroithms},(2005), 341-351.

\bibitem
{vinh3} Vinh, L. A., Random walk on hypergroup of conics over finite fields, to appear on
\textit{The Global Journal in Pure and Applied Mathematics}.

\bibitem
{quadrance} Wildberger, N. J.,
\textit{Divine Proportions: Rational trigonometry to universal geometry}, WildEgg, Australia, 2005.
\end{thebibliography}
\end{document}